\documentclass{pspum-l}

\newtheorem{theorem}{Theorem}[section]
\newtheorem{lemma}[theorem]{Lemma}
\usepackage[all]{xy}

\usepackage{amsmath, graphicx, tensor, stmaryrd, wasysym, amssymb}
\usepackage{amsfonts}
\usepackage{amssymb}
\usepackage{amscd}
\usepackage{amsthm}
\usepackage{latexsym}
\usepackage{amsbsy}
\usepackage{xcolor,enumerate}
\usepackage[normalem]{ulem}

\theoremstyle{definition}
\newtheorem{definition}[theorem]{Definition}

\theoremstyle{remark}
\newtheorem{remark}[theorem]{Remark}

\numberwithin{equation}{section}

\newcommand{\ot}{\otimes}

\newcommand{\Hc}{\mathcal{H}}
\newcommand{\Cc}{\mathcal{C}}
\newcommand{\Tc}{\mathcal{T}}
\renewcommand{\Mc}{\mathcal{M}}
\newcommand{\Zc}{\mathcal{Z}}

\newcommand{\Oc}{\mathcal{O}}
\newcommand{\Ac}{\mathcal{A}}

\newcommand{\hmod}{{_H\mathcal{M}}}

\newcommand{\hhmod}{{_{H^2}\mathcal{M}}}

\newcommand{\Dc}{\mathcal{D}}
\newcommand{\Nc}{\mathcal{N}}

\newcommand{\la}{\triangleright}
\newcommand{\ra}{\triangleleft}

\newcommand{\Bc}{\mathcal{B}}

\newcommand{\vect}{\text{Vec}}

\begin{document}

\title{Hopf Cyclic Cohomology and Beyond}

\author{Masoud Khalkhali}
\address{Department of Mathematics, University of Western Ontario, London, Ontario,
Canada, N6A 5B7}

\email{masoud@uwo.ca}
\thanks{The first author was supported in part by an NSERC Grant.}

\author{Ilya Shapiro}
\address{Department of Mathetics and Statistics,
University of Windsor, Windsor, Ontario, Canada, N9B 3P4}
\email{ishapiro@uwindsor.ca}
\thanks{The second author was supported in part by NSERC Grant \#820455.}


\dedicatory{}


\begin{abstract}
This paper is an introduction to Hopf cyclic cohomology with an emphasis on its most recent developments. We cover three major areas:  the original definition of Hopf cyclic cohomology by Connes and Moscovici as an outgrowth of their study of transverse index theory on foliated manifolds, the introduction of Hopf cyclic cohomology with coefficients by Hajac-Khalkhali-Rangipour-Sommerh\"auser,  and  finally the latest episode  on unifying the coefficients as well as extending the notion to more general settings beyond Hopf algebras.  In particular, the last section discusses the relative Hopf cyclic theory that arises in the braided monoidal category  settings. 
\end{abstract}

\maketitle


\tableofcontents

\section{Introduction}
Cyclic homology can be understood as a noncommutative analogue of de Rham cohomology of smooth manifolds. In the same vein,  Hopf cyclic cohomology  should be understood as  the noncommutative analogue of group homology and  Lie algebra homology. 
Having said so, it must be understood that there  was  no ``royal road"  to either theories and their discovery  came indirectly through  questions related to index theory of elliptic operators and its abstract formulation called  K-homology theory. Once cyclic cohomology was discovered, it  turned out that it is  possible to define  it in a purely algebraic setting  and in fact several  different, but equivalent,  variations are possible. The  original approach by Connes, based on quantized calculus,  and his route to the discovery of cyclic cohomology  is well documented in his 1981 Oberwolfach talk and his written report,  where cyclic cohomology was first announced \cite{ac81} (for quantized calculus see \cite{c2}):

\medskip

{\it ``The transverse elliptic theory for foliations requires as a
preliminary step a purely algebraic work, of computing for a
noncommutative algebra $\mathcal{A}$ the cohomology of the following
complex: $n$-cochains are multilinear functions \\$\varphi (f^0, \dots,
f^n)$ of $f^0, \dots ,f^n \in  \mathcal{A}$ where
$$\varphi (f^1, \dots,
f^0)=(-1)^n \varphi (f^0, \dots,
f^n)$$
and the boundary is
\begin{eqnarray*}
b\varphi (f^0, \dots,f^{n+1}) &=&\varphi (f^0f^1, \dots,f^{n+1})-\varphi (f^0, f^1f^2,\dots,
f^{n+1})+\dots \\& &+(-1)^{n+1} \varphi (f^{n+1}f^0, \dots,
f^n).
\end{eqnarray*}
The basic class associated to a transversally elliptic operator,
for $\mathcal{A}=$ the algebra of the foliation,  is given by:
$$ \varphi (f^0, \dots,f^n)=Trace \, (\varepsilon F[F, f^0][F, f^1]\dots[F, f^n]), \quad
f^i\in \mathcal{A}$$
where
$$F=\left(
\begin{matrix}
0& Q \\
P & 0
\end{matrix}
\right), \quad  \varepsilon= \left(
\begin{matrix}
1& 0 \\
0& -1
\end{matrix}
\right),
$$
and $Q$ is a parametrix of $P$. An operation
$$S: H^n(\mathcal{A})\to
H^{n+2}(\mathcal{A})$$
is constructed as well as a pairing
$$K(\mathcal{A})
\times H(\mathcal{A})\to \mathbb{C}$$ where $K(\mathcal{A})$ is
the algebraic $K$-theory of $A$. It gives the index of the
operator from its associated class $\varphi$. Moreover $\langle e, \,
\varphi \rangle = \langle e, \, S \varphi \rangle$ so that the important group to determine
is the inductive limit $H_p=\underset{\to}{Lim} H^n(\mathcal{A})$
for the map $S$. Using the tools of homological algebra the groups
$H^n(\mathcal{A}, \mathcal{A}^*)$ of Hochschild cohomology with
coefficients in the bimodule $\mathcal{A}^*$ are easier to
determine and the solution of the problem is obtained
in two steps:\\
1) the construction of a map
$$B: H^n(\mathcal{A},
\mathcal{A}^*) \to H^{n-1}(\mathcal{A})$$
and the proof of a long exact
sequence
$$
\cdots\to H^n(\mathcal{A}, \mathcal{A}^*) \overset{B}{\to} H^{n-1}(\mathcal{A})\overset{S}{\to}
H^{n+1}(\mathcal{A})\overset{I}{\to}H^{n+1}(\mathcal{A},
\mathcal{A}^*)\to \cdots
$$
where $I$ is the obvious map from the cohomology of the above
complex to
the Hochschild cohomology.\\
2) The construction of a spectral sequence with $E_2$ term given
by the cohomology of the degree $-1$ differential $I\circ B$ on
the Hochschild groups $H^n(\mathcal{A}, \mathcal{A}^*)$ and which
converges strongly to a graded group associated to the inductive
limit.

This purely algebraic theory is then used. For
$\mathcal{A}=C^{\infty}(V)$ one gets the de Rham homology of
currents, and for the pseudo-torus, i.e. the algebra of the
Kronecker foliation, one finds that the Hochschild cohomology
depends on the Diophantine nature of the rotation number while the
above theory gives $H^0_p$ of dimension $2$ and $H^1_p$ of dimension
$2$,  as expected, but from some
remarkable
cancelations."}

\medskip

Our goal in this paper is to give  a short survey of aspects of Hopf cyclic cohomology where we have been personally active in.    Hopf cyclic cohomology was defined by Connes and Moscovici in their seminal article
\cite{cm2} (see also \cite{cm3,cm4}) which  deals with index theory for transversally elliptic operators. The  resulting theory can be regarded as a noncommutative analogue of construction of characteristic classes in Chern-Weil theory through  group cohomology.   The role of group cohomology is now played by Hopf cyclic cohomology. Due to noncommutativity and non-cocommutativity of the Hopf algebras involved, defining this cohomology theory was not straightforward at all.    In fact one of the main steps in \cite{cm2} was
to obtain a
{\it noncommutative characteristic map}
 $$\chi_{\tau}: {HC}_{(\delta,\sigma)}^\bullet (H) \longrightarrow  HC^\bullet (A),$$
 for  an action of a  Hopf
algebra $H$ on an  algebra $A$  endowed with an  invariant trace  $\tau:
A \to \mathbb{C}$. Here, the pair $(\delta, \sigma)$ consists of a grouplike element $\sigma \in H$ and
a character $\delta: H\to \mathbb{C}$ satisfying certain compatibility conditions to be explained later in this paper. This  intriguing index theoretic  construction naturally begged   for its purely algebraic underpinnings to be abstracted  in order  to build  the theory independently of its  original index theory connections.

In \cite{kr2} a new approach was indeed  found by defining 
a cyclic cohomology theory  for triples $(C, H, M)$, where $C$ is  a coalgebra endowed with an action of
a Hopf algebra $H$ and  $M$ is an $H$-module and an $H$-comodule  satisfying
some extra compatibility conditions.
It was observed that
the theory  of Connes and Moscovici corresponds to $C=H$ equipped with
the regular action of $H$ and $M$ a one dimensional $H$-module with an extra
structure.

 One of the main ideas of \cite{kr2} was to view the
Hopf-cyclic cohomology as the  cohomology of the {\it invariant} part of
certain natural complexes attached
to $(C, H, M)$. This  should remind one of  the cohomology
 of the Lie algebra of a
 Lie group as the invariant part of the de Rham cohomology of the  group. Another important idea in \cite{kr2}   was to introduce {\it coefficients} into the theory.   This also explained the
important role played by the
so called modular pair in involution  $(\delta, \sigma)$ in \cite{cm2}.

The module $M$ is a noncommutative analogue of coefficients
 for  Lie algebra cohomology and group cohomology. And a natural question  was to identify   the most general type of modules. Now  the periodicity
condition $\tau_n^{n+1}=id$ for the cyclic operator and the fact that all simplicial and
cyclic operators have to descend to the invariant complexes, puts restrictions 
on the type of the $H$-modules $M$ one can work with. This problem
that was partly 
 solved  in   \cite{kr2}   was then  completely solved
 in the twin papers of  Hajac-Khalkhali-Rangipour-Sommerh\"auser  \cite{hkrs1, hkrs2} by
 introducing the class of {\it stable anti-Yetter-Drinfeld modules} over a
 Hopf algebra. It was shown that the category of anti-Yetter-Drinfeld modules over a Hopf algebra $H$ is obtained from the category of  Yetter-Drinfeld $H$-modules by twisting the latter by a modular pair in involution. 
 

Starting with Section \ref{nonsense} we focus on developing a conceptual understanding of the coefficients, i.e., (stable) anti-Yetter-Drinfeld modules and their generalizations.  The main benefit of this approach, besides a general method of extending the definition to other Hopf-like or bialgebra-like settings, is that it suggests that stable anti-Yetter-Drinfeld modules, or rather their  mixed generalizations, are exactly the $D$-modules on the noncommutative ``space" obtained as a quotient of a point by the symmetries given by a Hopf algebra.   We see this by observing that stable anti-Yetter-Drinfeld modules arise as the naive cyclic homology category of the category of representations of a Hopf algebra.  On the other hand \cite{bfn} shows that $HH(QCohX)^{S^1}$, the cyclic homology category, in the $\infty$-setting, of the category of quasi-coherent sheaves is, roughly speaking, the category of $D$-modules.  

Once the correct abelian category of coefficients is found, and the algebras are mapped to it, the cohomology is obtained as an $Ext$, just as in \cite{ac83b}.  More precisely, given an $A\in\Cc$ a unital associative algebra, and $M\in HH(\Cc)^{S^1}$ a coefficient, the definition of the analogue of Hopf-cyclic cohomology of $A$ with coefficients in $M$ can be obtained as \cite{shcatchern}: $$HC^\bullet(A,M)=Ext^\bullet_{HH(\Cc)^{S^1}}(ch(A),M),$$ where $ch(A)$ is defined in \eqref{cherndefi}.  The analogue of the deRham cohomology of $M$ would then be $Ext^\bullet_{HH(\Cc)^{S^1}}(ch(1),M)$. We will not actually deal with such cohomology groups in this paper outside of Theorem \ref{esl}.  Our focus instead will be on understanding the coefficients, namely the categories $HH(\Cc)$, $HC(\Cc)$, and a bit of $HH(\Cc)^{S^1}$; the latter does not appear outside of Sections \ref{sigmacyclic} and \ref{hmor}.   We note that the naive treatment of Section \ref{nonsense} can be moved to the $\infty$-setting for a more direct comparison with \cite{bfn}.

All  algebras, coalgebras and Hopf algebras in this paper  are over a fixed  field $k$ of
characteristic zero, and are unital or counital. Tensors and homs are denoted   $\otimes,  Hom$  and will be over k, unless specified otherwise.  
We use  Sweedler's  notation for
comultiplication $\Delta$, with summation understood, and write
$$ \Delta (c)= c^{(1)}\otimes c^{(2)}.$$
The   {\it iterated comultiplication} maps
$$\Delta ^n =(\Delta \otimes I)\circ \Delta^{n-1}: C\longrightarrow C^{\otimes (n+1)},\quad \Delta^1=\Delta,$$
will be written as 
$$\Delta^n (c)= c^{(1)}\otimes \cdots \otimes c^{(n+1)},$$
where summation is understood.

\section{Hopf algebras in cyclic homology}

Even before the introduction of Hopf cyclic cohomology, Hopf algebras  played a role in  cyclic homology theory, acting as symmetries of  associative algebras.  In this section we shall briefly look at two such examples. 

\subsection{Deformation complex and Deligne's conjecture}
In 1960's Gerstenhaber understood that deformations of an associative algebra $A$ is controlled 
by   Hochschild cohomology of $A$ with coefficients in $A$, usually denoted as $H^\bullet (A, A)$ \cite{Ger}. In particular he showed that $H^\bullet (A, A)$ is a graded Poisson algebra, i.e.,  
a commutative differential graded algebra   equipped with a compatible differential graded Lie algebra structure    with a shifted  degree. A structure that is usually called a Gerstenhaber algebra. He had a similar result for deformations of Lie algebras. The product and the Lie bracket operations are defined  on the level of cochains and we have a cochain complex   $(C^\bullet (A, A), \delta, \cup,  [ \, , \, ])$, where 
  $$C^n (A, A) = Hom_k (A^{\otimes n}, A).$$ 
      The Jacobi identity holds  for cochains, but the cup product is only homotopy associative.  Compatibility of the cup product with Lie algebra structure  is also only up to  homotopy. This structure can be rather loosely called a   homotopy Gerstenhaber algebra. And a major question   is  what is the full algebraic structure hidden here?
 Deligne's conjecture gives a precise answer to this question:   The Hochschild  cochain complex  is an  algebra  over the
singular chain operad of the little squares operad $E_2$ \cite{Del}.

In \cite{kh1} an algebraic approach to this question is developed. It uses  a  differential graded Hopf algebra  that  acts as symmetries of the bar complex and a bar-cobar duality.  A central tool here
is the notion of $X-$complex for differential graded  coalgebras  as in \cite{Q} (see also [3] for the case of algebras). Let us explain this construction briefly.  Recall that the bar construction is a functor from the category of differential graded algebras to the category of differential graded coalgebras
$$\mathcal{B}: \mathcal{DGA} \to \mathcal{DGC},$$ 
where  $ \mathcal{B} (A) = \oplus_n A^{\otimes n}$. 
A crucial property of the bar constructions is that its coderivations is isomorphic to the deformation complex:
 $$\text{Coder} (\mathcal{B} (A), \mathcal{B} (A)) = C^\bullet (A, A).$$
It is shown in \cite{kh1} that   If $A$ is a homotopy Gerstenhaber algebra, then $\mathcal{B} (A)$ is a Hopf algebra.This is how Hopf algebras emerge as relevant to operations on cyclic homology. 

The $X$-complex of a differential graded coalgebra $C$, denoted $\hat{X} (C)$ is the $\mathbb{Z}/2\mathbb{Z}$-graded complex 
$$ C \to \Omega^1 C_{\flat}$$
where $\Omega^1 C_{\flat}$ is the cocommutator subspace of the $C$-bicomodule of universal differential 1-forms on $C$. 
By a   result
of Quillen \cite{Q}, the cyclic  chain complex of an algebra $A$   is isomorphic to  the $X$-complex of the
bar construction
$$   \mathcal{B} (A) \to   \Omega^1 (\mathcal{B} (A))_{\flat}.$$
From this point of view operations on cyclic and Hochschild complex of
homotopy Gerstenhaber  algebras are predicted by Kunneth formula for the $X$-complex of differential
graded coalgebras.

More precisely, let $V$ be a homotopy Gerstenhaber algebra. Then it is shown in \cite{kh1} that there are natural maps of
supercomplexes
$$ \hat{X}(BV)\otimes \hat{X}(BV) \longrightarrow \hat{X}(BV),\quad \hat{X}(BV)\otimes \hat{X}(BV_0)\longrightarrow \hat{X}(BV_0).$$  Applied to $V= C^\bullet (A, A)$, one can get many of the operations.

   The main point here is that Hopf algebra symmetries lead to operations in cyclic homology and to a host of intriguing relations among them. We should also mention that the bar construction and Quillen's theorem above should be useful to define cyclic homology of $A_{\infty}$-algebras and algebras over operads in general. The point is that  the bar construction of an  $A_{\infty}$-algebra is a differential graded coalgebra and hence one can consider its $X$-complex as its cyclic complex. For a recent survey of the status of   Deligne's conjecture the reader can consult the paper of Kaufmann and Zhang \cite{kz}.

\subsection{ Cyclic homology of Hopf  crossed products}

An  important question in cyclic cohomology theory is to compute the cyclic
cohomology of a crossed product algebra  $ A\sharp G$  where the group $G$ acts
on the algebra $A$ by automorphisms. This question  has been studied by many  since the beginning of the theory in 1980's.  If $G$  is a discrete group, there is a spectral sequence, due to
Feigin and Tsygan \cite{FT}, which converges to the cyclic homology of the crossed
product algebra. This result generalizes Burghelea’s calculation of the cyclic
homology of a group algebra \cite{bur}.
For a recent survey we refer to \cite{ponge} in this volume and references therein.

 In \cite{ak1} a spectral sequence is obtained  which computes the cyclic cohomology of a Hopf crossed product algebra. The Hopf algebra need not be commutative or cocommutative. This result was later extended to compute the Hopf cyclic cohomology of bicrossed products of two Hopf algebras in \cite{moscrang1, moscrang2}. This gave a uniform method to compute the Hopf cyclic cohomology of Connes-Moscovici Hopf algebras in all dimensions since the latter  are known to be Hopf bicrossed product algebras \cite{cm2}.

 Let a Hopf algebra $\mathcal{H}$  act on an algebra $A$ such that the coalgebra structure of $H$ is compatible with the algebra structure of $A$ and  $A$ is an  $\mathcal{H}$-module. We say $A$ is a Hopf module algebra. We can then form the Hopf crossed product algebra $A \sharp \mathcal{H}$.  In \cite{ak1} it is shown that   there exists an isomorphism
between $ {C}_{\bullet}(A^{op} \sharp  \mathcal{H}^{cop})$ the cyclic module of the crossed product algebra $A^{op} \sharp  \mathcal{H}^{cop} $, and $\Delta(A \natural \mathcal{H}) $, the cyclic
module defined as  the diagonal of a  so called cylindrical module  $A \natural \mathcal{H}$.
When the antipode of the Hopf algebra  is invertible, it was shown that the cyclic homology  $HC_{\bullet}(A \sharp  \mathcal{H})$ can be approximated by a spectral
sequence and     the of $ {E}^0 , {E}^1$ and $ {E}^2 $  of  this
spectral sequence was explicitly computed.

\section{The work of Connes and Moscovici}

In this section we recall the approach by  Connes and Moscovici
towards their  definition of a 
cyclic cohomology theory for Hopf algebras. 
The local
index formula \cite{cm1} gives the Connes-Chern character
$Ch (A, h, D)$ of a regular
spectral triple $(A, h, D)$ as a  cyclic cocycle in the $(b, B)$-bicomplex of the algebra $A$.
For  spectral
triples of interest in transverse geometry \cite{cm2}, this cocycle is in  the
image of the Connes-Moscovici characteristic map $\chi_{\tau}$ defined below.   To identify
this class in  terms of characteristic classes of foliations, it turned out that it would be
 helpful to
 show that it is the image of a cocycle for a cohomology theory for Hopf algebras. This is rather
similar to the situation for classical characteristic classes which are pull backs of group
cohomology classes via a classifying map. 

We can formulate this problem abstractly as follows. Let $H$ be a Hopf algebra, $\delta: H\to k$ a character and
$\sigma \in H$ a grouplike element. Following \cite{cm2, cm3, cm4}, we
say that $(\delta, \sigma)$ is a {\it modular pair} if $\delta
(\sigma)=1$, and a {\it modular pair in involution} if 
$$ \widetilde{S}_{\delta}^2 (h)=\sigma
h\sigma^{-1},$$
for all $h$ in $H$. Here the
  $\delta$-{\it twisted antipode}
$\widetilde{S}_{\delta}: H \rightarrow H$ is defined as 
$$\widetilde{S}_{\delta}(h)= \delta (h^{(1)})S(h^{(2)}).$$

Now let  $A$ be an $H$-module algebra and  $\tau: A \to
k$ be a $\delta$- invariant $\sigma$-trace on $A.$ 
 The Connes-Moscovici {
characteristic map is defined as 
$$\chi_{\tau}: H^{\otimes n}\longrightarrow Hom (A^{\otimes (n+1)},\;
k),$$
\begin{eqnarray*}
\chi_{\tau}(h_1\otimes \cdots \otimes h_n)(a_0\otimes \cdots \otimes a_n)=
\tau (a_0 h_1(a_1)\cdots h_n(a_n)).
\end{eqnarray*}
This is a generalization of a map defined  for an action of a Lie algebra  on an algebra in \cite{ac80} which was then fully developed in \cite{c2}.  The question  is if one can  define a cocyclic module structure on  the collection of
spaces $\{H^{\otimes n}\}_{n\geq 0}$ 
such that the characteristic map $\chi_{\tau}$ turns into a
morphism of cocyclic modules? The face, degeneracy,
and cyclic operators for  $Hom (A^{\otimes (n+1)},\; k) $ are
defined as 
\begin{eqnarray*}
\delta^n_i \varphi (a_0, \cdots, a_{n+1})&=&\varphi (a_0, \cdots,
a_ia_{i+1}, \cdots, a_{n+1}), \quad i=0, \cdots, n,\\
\delta^n_{n+1} \varphi (a_0, \cdots, a_{n+1})&=&\varphi (a_{n+1}a_0,
a_1, \cdots, a_n),\\
\sigma^n_i\varphi (a_0, \cdots, a_{n}) &=&\varphi (a_0, \cdots,a_i, 1,
\cdots, a_{n}),\quad i=0, \cdots, n,\\
\tau_n \varphi (a_0, \cdots, a_n)& =&\varphi (a_n, a_0, \cdots, a_{n-1}).
\end{eqnarray*}
The
relation
$$h(ab)=h^{(1)}(a)h^{(2)}(b)$$
shows that  the face operators
on $H^{\otimes n}$ must involve the coproduct of
$H$: 
\begin{eqnarray*}
\delta^n_0(h_1\otimes
\cdots \otimes h_n)&=& 1\otimes h_1\otimes
\cdots \otimes h_n,\\
 \delta^n_i(h_1\otimes
\cdots \otimes  h_n)&=& h_1\otimes \cdots \otimes h_i^{(1)}\otimes h_i^{(2)}\otimes \cdots
\otimes h_n,\\
\delta^n_{n+1}(h_1\otimes \cdots \otimes h_n)&=&h_1\otimes
\cdots \otimes  h_n\otimes \sigma,
\end{eqnarray*}
With this definition, we have 
$$\chi_{\tau}\delta_i^n =\delta_i^n \chi_{\tau}.$$
Similarly, the degeneracy operators on $H^{\otimes n}$ should
be  defined as 
$$\sigma^n_i(h_1\otimes \cdots \otimes h_n)=h_1 \otimes \cdots \otimes \varepsilon
(h_i)\otimes \cdots \otimes h_n.$$
A  challenging part is to define the 
 {\it cyclic operator} $\tau_n : H^{\otimes n} \to H^{\otimes n}$.
 Compatibility of the latter  with $\chi_{\tau}$ implies the relation 
$$ \tau (a_0  \tau_n(h_1\otimes \cdots \otimes h_n)(a_1\otimes \cdots \otimes a_n))=
\tau (a_n  h_1(a_0)h_2(a_1)\cdots h_n(a_{n-1})).$$
Now  the  $(\delta, \sigma)$-invariance property of  $\tau$  shows that 
$$\tau (a_1 h(a_0))=\tau (h(a_0)\sigma (a_1))=\tau (a_0\tilde{S}_{\delta}(h)(\sigma (a_1)).$$
This suggests a definition for  $\tau_1: H \to H$ as
$$\tau_1(h)= \tilde{S}_{\delta}(h)\sigma.$$
Note that the condition
$\tau_1^2=I$ is equivalent to the involutive condition
$\tilde{S}_{\delta}^2 (h) =\sigma  h \sigma^{-1}.$
In general one can show that for any $n$:
\begin{eqnarray*}
\tau (a_n h_1(a_0) \cdots h_n(a_{n-1}))& =&\tau (
h_1(a_0) \cdots h_n(a_{n-1})\sigma (a_n))\\
&=&\tau(a_0 \tilde{S}_{\delta}(h_1)(h_2 (a_1)\cdots h_n(a_{n-1})\sigma (a_n)))\\
&=& \tau(a_0 \tilde{S}_{\delta}(h_1)\cdot(h_2 \otimes \cdots \otimes h_n\otimes \sigma)
(a_1\otimes \cdots \otimes a_n).
\end{eqnarray*}
This suggests that the {\it Hopf-cyclic operator} $\tau_n : H^{\otimes
n}\to  H^{\otimes n}$ should be defined as
$$\tau_n (h_1\otimes \cdots \otimes h_n)=\tilde{S}_{\delta}(h_1)\cdot(h_2 \otimes \cdots \otimes
h_n\otimes \sigma),$$
where $\cdot$ denotes the diagonal action defined by
$$h\cdot (h_1\otimes \cdots \otimes h_n):= h^{(1)}h_1\otimes
h^{(2)}h_2\otimes \cdots \otimes h^{(n)}h_n.$$ We let $\tau_0=I:
H^{\otimes 0}=k \to H^{\otimes 0}$, be the identity map. 

 Connes and Moscovici proved the remarkable result  that 
the above face, degeneracy, and cyclic
operators on  $\{H^{\otimes n} \}_{n\geq 0}$ define a cocyclic module  \cite{cm2, cm3, cm4}. 
The resulting cyclic cohomology groups are denoted by
$HC^n_{(\delta, \sigma)}(H)$, $n=0, 1,\cdots$ and we have  the
desired  characteristic map
$$\chi_{\tau}: HC^n_{(\delta, \sigma)}(H) \to HC^n(A).$$

The Hopf cyclic cohomology groups $HC^n_{(\delta, \sigma)}(H)$ are
computed for several key examples in \cite{cm2, cm3, cm4, moscrang1, moscrang2}.  For a recent survey and references we refer to the article of Moscovici in this volume \cite{mos}. In  particular 
for applications to transverse index theory,  
 the (periodic)  cyclic cohomology of the Connes-Moscovici Hopf
algebra $\mathcal{H}_1$ plays an important role.  It is shown in \cite{cm2} that the
periodic groups $HP^n_{(\delta, 1)}(\mathcal{H}_1)$ are canonically
isomorphic to the Gelfand-Fuchs cohomology of the Lie algebra of
formal vector fields $\mathfrak{a}_1$ on the line:
$$H^\bullet(\mathfrak{a}_1, \mathbb{C})=HP^\bullet_{(\delta,
1)}(\mathcal{H}_1).$$
The following
interesting elements appear. 
It can be directly checked that the elements $\delta_2' :
=\delta_2-\frac{1}{2}\delta_1^2$ and $\delta_1$ are cyclic 1-cocycles
on $\mathcal{H}_1$, and
$$F:=X\otimes Y-Y\otimes X -\delta_1Y\otimes Y$$
is a  cyclic 2-cocycle. See \cite{como6} for a survey  and  relations
between these cocycles and
 the Schwarzian derivative, Godbillon-Vey cocycle, and
 the transverse fundamental class of Connes \cite{c3}, respectively.

\section{Hopf cyclic cohomology with coefficients}

In \cite{kr2} (see also \cite{ak2} for an equivariant version) the Hopf cyclic cohomology was extended to a so called  {\it invariant  cyclic cohomology} theory for triples $(C, H, M)$ where $C$ is a $H$-module coalgebra and  $M$ is an $H$-module. The three pieces $C, H, M$  had to satisfy some compatibility conditions so that the resulting complex turns out to be a cyclic module. An important question which was left open in   \cite{kr2} was to identify 
     the most general class of coefficients  $M$ for  this Hopf cyclic cohomology. This
    problem
    was   completely solved   in the twin papers of Hajac-Khalkhali-Rangipour-Sommerhauser  \cite{hkrs1, hkrs2}. It was  shown  that
    the most general coefficients are the
    class of  so called {\it  stable anti-Yetter-Drinfeld modules} (AYD modules). It was found that anti-Yetter-Drinfeld modules are twistings by modular pairs in
involution  of Yetter-Drinfeld  modules.  This means that the category of
anti-Yetter-Drinfeld modules is a `mirror' of the category
of Yetter-Drinfeld  modules. 
It is quite surprising that  when the general formalism of cyclic
cohomology theory, namely the theory of (co)cyclic modules
\cite{ac83b},  is applied to Hopf algebras, variations of  
 notions like Yetter-Drinfeld  modules which are of importance in low dimensional topology and  invariants of knots and ribbons appear
naturally.

\subsection{Stable anti Yetter-Drinfeld modules}

Let $H$ be a Hopf algebra,  $M$  a left $H$-module and left $H$-comodule. Then 
 $M$ is  called a left-left {\it Yetter-Drinfeld $H$-module} if 
$$\rho(hm)=h^{(1)}m^{(-1)}S(h^{(3)})\ot h^{(2)}m^{(0)},$$
 for all $h\in H$ and $m\in M$ \cite{maj, rt, y}.  We denote the category of left-left Yetter-Drinfeld  modules over $H$ by
 $^{H}_{H} \mathcal{Y}\mathcal{D}$.
This notion is closely related to the
 Drinfeld double of finite dimensional Hopf algebras. In fact if $H$ is
 finite dimensional, then one  can show that the category $^{H}_{H} \mathcal{Y}\mathcal{D}$  is
 isomorphic to the category of left modules over the Drinfeld double
 $D(H)$ of $H$ \cite{maj}.

 The category $^{H}_{H} \mathcal{Y}\mathcal{D}$ is a braided monoidal category under the
braiding
$$c_{M, N}: M\otimes N \to N\otimes M, \quad m\otimes n \mapsto m^{(-1)}\cdot n\otimes  m^{(0)},$$   
    and  the  tensor product $M\otimes N$ of two Yetter-Drinfeld  modules defined by structure maps 
$$h (m\otimes n)=h^{(1)}m\otimes h^{(1)}n, \quad (m\otimes n)\mapsto
m^{(-1)}n^{(-1)}\otimes m^{(0)}\otimes n^{(0)}.$$

 The category $^{H}_{H} \mathcal{Y}\mathcal{D}$ is the {\it
center} of the monoidal category  $\hmod$ of left $H$-modules. Recall that the (left)
center $\mathcal{Z}\mathcal{C}$ of a monoidal category \cite{kas} is
a category whose objects are pairs $(X, \sigma_{X, -})$, where $X$
is an object of $\mathcal{C}$ and $\sigma_{X, -}: X\otimes -\to
-\otimes X$ is a natural isomorphism satisfying certain
compatibility axioms with associativity and unit constraints. It
can be shown that the center of a monoidal category is a braided
monoidal category and
$\mathcal{Z} (\hmod)=^{H}_{H} \mathcal{Y}\mathcal{D}$  \cite{kas}.

 Let us give an example of a class of Yetter-Drinfeld modules. Let $H=kG$ be the group algebra of a discrete group $G$.  Let $\mathcal{G}$ be a groupoid whose objects are $G$ and its morphisms are defined by
$$Hom(g, h)=\{k\in G; \; kgk^{-1}=h\}.$$
 Recall that an {\it action}
 of a groupoid $\mathcal{G}$ on the category $\vect$ of vector spaces is simply a
 functor $F: \mathcal{G} \to \vect$.
  Then it is easily seen that we have a
one-one correspondence between YD modules for $kG$ and actions of $\mathcal{G}$ on $\vect$. This example
clearly shows that one can think of an YD module over $kG$ as an {\it equivariant sheaf } on $G$ and
of $YD$ modules as noncommutative analogues of equivariant sheaves on a topological  group.  See Section \ref{Hopfex} for the dual case of $H=\Oc_G$, i.e, functions on a monoid $G$.

 Unlike Yetter-Drinfeld modules, the definition of anti-Yetter-Drinfeld modules  and stable anti-Yetter-Drinfeld modules was entirely motivated   by cyclic
cohomology theory. In fact the  anti-Yetter-Drinfeld condition guarantees
that the simplicial and cyclic operators are well defined on
invariant  complexes and the stability condition implies
that the crucial periodicity condition for cyclic modules are
satisfied.  A center-like description is possible, see Section \ref{sec:bi}.

\begin{definition}\label{st11} (\cite{hkrs1}) A left-left  anti-Yetter-Drinfeld  $H$-module is a left $H$-module and
left $H$-comodule such that
$$\rho(hm)=h^{(1)}m^{(-1)}S^{-1}(h^{(3)})\ot h^{(2)}m^{(0)},$$
for all $h\in H$ and  $m\in M.$
 $M$ is called stable if in addition we have
$$m^{(-1)}m^{(0)}=m,  $$
for all $m\in M$.
\end{definition}

The following lemma from \cite{hkrs1} shows that 1-dimensional
SAYD modules correspond to Connes-Moscovici's modular pairs in
involution:
\begin{lemma} There is a one-one correspondence between modular pairs in
involution $(\delta, \sigma)$ on $H$ and SAYD module structure on $M=k$, defined by
$$h. r=\delta (h) r, \quad r\mapsto \sigma \otimes r,$$
for all $h\in H$ and $r\in k$. We denote this module by $M=^{\sigma}\!\!\!k_{\delta}.$
\end{lemma}

Let $^{H}_{H} \mathcal{A}\mathcal{Y}\mathcal{D}$ denote the
category of left-left anti-Yetter-Drinfeld $H$-modules, where
morphisms are $H$-linear and $H$-colinear maps. Unlike YD modules,
anti-Yetter-Drinfeld modules do not form a monoidal category under
the standard tensor product. This can be checked easily on
1-dimensional modules given by modular pairs in involution.
 The following result of  \cite{hkrs1}, however,
shows that the tensor
product of an anti-Yetter-Drinfeld module with a Yetter-Drinfeld module
is again anti-Yetter-Drinfeld.
\begin{lemma}\label{lem43}
Let $M$ be a Yetter-Drinfeld module and   $N$ be an anti-Yetter-Drinfeld module (both left-left).
  Then $M \otimes N $ is an anti-Yetter-Drinfeld module under the
  diagonal action and coaction:
  $$h (m\otimes n)=h^{(1)}m\otimes h^{(1)}n, \quad (m\otimes n)\mapsto
m^{(-1)}n^{(-1)}\otimes m^{(0)}\otimes n^{(0)}.$$
\end{lemma}

In particular, using a modular pair in involution $(\delta, \sigma)$, we
obtain a full and faithful functor
$$^{H}_{H} \mathcal{Y}\mathcal{D} \to ^{H}_{H} \mathcal{A}\mathcal{Y}\mathcal{D},
\quad M\mapsto \overset{-}{M}= ^{\sigma}\!\!\!k_{\delta}\otimes M.$$

This result clearly shows that AYD modules can be regarded as the
{\it twisted analogue}  or {\it mirror image} of YD modules, with
twistings provided by modular pairs in involution.

 Noncommutative
 principal bundles,  also known as Hopf-Galois extensions, 
 give rise to large
classes of examples of SAYD modules \cite{hkrs1}. Let $P$ be a right
$H$-comodule algebra, and let
$$B:=P^H=\{p\in P; \, \rho(p)=p\otimes 1\}$$
be the space of coinvariants of $P$. It is easy to see that $B$ is
a subalgebra of $P$. The extension $B\subset P$ is called a
Hopf-Galois extension if the map
$$can:  P\otimes_B P\to B\otimes H, \quad p\otimes p'\mapsto p\rho (p'), $$
is bijective. The bijectivity assumption allows us to
define the translation map $T: H \rightarrow P\ot_B P$,
$$T(h)=can^{-1}(1\ot h)=h^{(\bar{1})}\ot h^{(\bar{2})}.$$
It can be
checked that
 the centralizer  $Z_B(P)=\{p\mid bp=pb~~ \forall b\in B\}$ of $B$
 in $P$ is a subcomodule of $P$. There is an action of $H$ on
 $Z_B(P)$ defined by $ph=h^{({1})} p h^{({2})}$ called the
 Miyashita-Ulbrich action. It is shown that this action and
 coaction satisfy the Yetter-Drinfeld compatibility condition. On
 the other hand if $B$ is central, then by defining the new action $ph=(S^{-1}(h))^{({2})} p
 (S^{-1}(h))^{({1})}$and the right coaction of $P$ we have a SAYD
 module.

\subsection{Hopf cyclic cohomology with coefficients}

In the work of Hajac-Khalkhali-Rangipour-Sommerhaeuser \cite{hkrs2}  3 types of Hopf-cyclic cohomology theory are defined. In this section we give a brief account of two of them.

Let $A$ be a left $H$-module algebra  and  $M$
be a left-right SAYD  $H$-module.
 Then the spaces $M\otimes A^{\otimes (n+1)}$ are right $H$-modules
via the diagonal action
$$(m\otimes \widetilde{a})h:=mh^{(1)}\otimes S (h^{(2)})\widetilde{a},$$
where the left  $H$-action on $\widetilde{a} \in A^{\otimes (n+1)}$ is
 via the left diagonal action of $H$.

We  define the space of {\it equivariant  cochains on $A$ with
coefficients in $M$} by
$$\mathcal{C}^n_H(A, M):= Hom_H(M\otimes A^{\otimes (n+1)}, k).$$
More explicitly, $f: M\otimes A^{\otimes (n+1)} \rightarrow k$ is in
$\mathcal{C}^n_H(A, M)$, if and only if
$$f((m\otimes a_0\otimes \cdots \otimes a_n)h)= \varepsilon (h) f(m\otimes a_0\otimes \cdots \otimes a_n), $$
for all $h \in H, m\in M$, and $a_i\in A$. It is shown in \cite{hkrs2}
that the following operators define a cocyclic module structure on $\{\mathcal{C}^n_H(A, M)\}_{n\in \mathbb{N}}$:
\begin{align*}&
(\delta_if)(m\otimes a_0\otimes \cdots\otimes a_n)
= f(m\otimes a_0\otimes \cdots \otimes a_i a_{i+1}\otimes\cdots \otimes a_n),
~~~ 0 \leq i <n,~~~~~~\\ &
(\delta_nf)(m\otimes a_0\otimes \cdots\otimes a_n)
=f(m^{(0)}\otimes (S^{-1}(m^{(-1)})a_n)a_0\otimes a_1\otimes\cdots \otimes a_{n-1}),~~~~~~
\\ &
(\sigma_if)(m\otimes a_0\otimes \cdots\otimes a_n)
= f(m\otimes a_0 \otimes \cdots \otimes a_i\otimes 1 \otimes \cdots \otimes a_n),
~~~0\le i\le n,~~~~~~
\\  &
(\tau_nf)(m\otimes a_0\otimes \cdots\otimes a_n)
= f(m^{(0)}\otimes S^{-1}(m^{(-1)})a_n\otimes a_0\otimes \cdots \otimes a_{n-1}).~~~~~~
\end{align*}

We denote the resulting cyclic cohomology theory by $HC^n_H(A, M), n=0,
1, \cdots.$
For $M=H$ and $H$ acting on $M$ by conjugation and coacting
via  coproduct (Example 2.2.3.), we obtain the equivariant cyclic
cohomology theory of Akbarpour and Khalkhali for $H$-module algebras \cite{ak1, ak2}.

 Here is another variation.  
Connes-Moscovici's original example of Hopf-cyclic cohomology
belongs to this class of theories.
Let $C$ be a left $H$-module coalgebra, and $M$ be a right-left SAYD $H$-module. Let
$${\mathcal C}^n(C,M):=M\otimes_H C^{\otimes(n+1)} \quad n\in{\mathbb N}.$$
It can be checked that, thanks to the SAYD condition on $M$,   the
following operators are well defined and define a cocyclic module, denoted
$\{{\mathcal C}^n(C,M)\}_{n\in{\mathbb N}}$. In particular the
crucial periodicity conditions $\tau_n^{n+1}=id, \, n=0, 1, 2 \cdots$,
are  satisfied \cite{hkrs2}:
\begin{eqnarray*}
\delta_i(m\otimes c_0 \otimes \cdots \otimes c_{n-1})
&=&
m\otimes  c_0 \otimes\cdots\otimes  c_i^{(1)}\otimes c_i^{(2)}\otimes
c_{n-1},\,
 0 \leq i <n,\\
\delta_{n}(m\otimes c_0 \otimes \cdots \otimes c_{n-1})
&=&
m^{(0)}\otimes  c_0^{(2)}\otimes c_1
\otimes \cdots \otimes c_{n-1}  \otimes m^{({-1})}c_0^{(1)},
\\
\sigma_i(m\otimes c_0 \otimes \cdots \otimes c_{n+1})
&=&
m\otimes c_0 \otimes \cdots
\otimes \varepsilon(c_{i+1})\otimes
\cdots\otimes c_{n+1},\,
0\leq i \leq n,\\
\tau_n(m\otimes c_0 \otimes  \cdots \otimes c_n)
&=&
m^{(0)}\otimes c_1 \otimes
\cdots \otimes c_n \otimes m^{({-1})}c_0.
\end{eqnarray*}

 For $C=H$ and $M=^{\sigma}\!\!\!k_{\delta}$, the
cocyclic module $\{{\mathcal C}^n_H(C,M)\}_{n\in{\mathbb N}}$ is
isomorphic to  the cocyclic module of Connes-Moscovici \cite{cm2},
attached to a Hopf algebra endowed with a modular pair in
involution. This example is truly fundamental and started the
whole theory.

It is possible to interpret the stable anti-Yetter-Drinfeld $H$-modules as a local system over $H$, that is a module equipped with a noncommutative flat connection \cite{KayKha1}. Excision in Hopf cyclic cohomology is studied in \cite{KayKha2}. Other references include \cite{Gor, kh2, kr6, kr3, kr4, kr5}.

\section{Beyond Hopf algebras, a brief detour into limits}\label{nonsense}
Let $I$ be a small category (\cite{maclane} should suffice for this section).  Recall that an $I$-indexed diagram in a category $\Tc$ is a functor  $D_\bullet: I\to \Tc$.  These naturally form the category $Fun(I,\Tc)$.  We have an obvious functor $c:\Tc\to Fun(I,\Tc)$ that assigns to an object $T\in\Tc$ the constant diagram with $c(T)_i=T$ for all $i\in I$ and for any $f\in Hom_I(i,j)$, the induced map $f_*\in Hom_\Tc(c(T)_i,c(T)_j)$ is $Id_T$.  Should they exist, we denote by $\varinjlim_I$ the left adjoint of $c$ and by  $\varprojlim_I$ the right adjoint of $c$ to denote the colimit and limit functors; the use of the latter will be delayed until Section \ref{ulc}.

Since we are interested in the $I$-indexed diagrams in the $2$-category $\Cc at$ of small categories enriched over $\vect$, i.e.,  vector spaces, let us briefly spell out the subtleties that this categorification entails.  A diagram $\Cc_\bullet:I\to \Cc at$ consists of categories $\Cc_i$, an assignment of $f_*\in Fun(\Cc_i,\Cc_j)$ to every $f\in Hom_I(i,j)$, and natural isomorphisms $\iota_{f,g}\in NatIso(f_*g_*,(fg)_*)$.  The latter satisfy $$\iota_{f,gh}\iota_{g,h}=\iota_{fg,h}\iota_{f,g}.$$  We can define the weak version of the above by relaxing the condition that $\iota_{f,g}$'s be isomorphisms.  

Given two diagrams $\Cc_\bullet$ and $\Dc_\bullet$ a morphism of diagrams is the data $(F_\bullet,\beta_f)$, where $F_i\in Fun(\Cc_i,\Dc_i)$ and $\beta_f\in NatIso(f^\Dc_*F_i, F_jf^\Cc_*)$ for $f\in Hom_I(i,j)$.  The $\beta_f$'s satisfy $$\iota^\Cc_{f,g}\beta_f\beta_g=\beta_{fg}\iota^\Dc_{f,g}.$$ A weak morphism is obtained by again relaxing the condition that $\beta_f$'s be isomorphisms.  

A constant diagram $c(\Tc)_\bullet$ is obtained from a category $\Tc$ in an obvious manner; namely,  $c(\Tc)_i=\Tc$, $f_*=Id_\Tc$, and $\iota_{f,g}=Id$. By a (weak) section of $\Cc_\bullet$ we mean a (weak) morphism from $c(\vect)_\bullet$ to $\Cc$.  Explicitly, it is a choice of $C_i\in \Cc_i$ together with $\beta_f\in Hom_{\Cc_j}(f_*C_i,C_j)$, where the $\beta_f$ is an isomorphism if the section is not weak.  The compatibility condition becomes: $$\beta_f f_*(\beta_g)=\beta_{fg}\iota_{f,g}.$$ 

It is immediate that (weak) morphisms map (weak) sections to (weak) sections, furthermore a weak section of $c(\Tc)_\bullet$ is just an object in $Fun(I,\Tc)$, while a strong section is much less reasonable, requiring all maps in $I$ to be mapped to isomorphisms in $\Tc$.  Thus, as tempting as it is to ignore the weak versions of the above concepts, one does need them as we shall see below.  Of course a strong version of any of the above may freely be used in place of a weak one.

Let us again define $\varinjlim_I$ to be the left adjoint to $c$, the inclusion of constant diagrams into (weak) diagrams.  According to the above discussion, we see that if our goal (and it is, to first order) is to obtain objects in $Fun(I,\Tc)$ from weak sections in $\Cc_\bullet$ we need only obtain an object in $Fun(\varinjlim_I\Cc_\bullet,\Tc)$.  In particular, given an $M\in \varinjlim_I\Cc_\bullet$ we obtain an object in $Fun(I^{op},\vect)$ via $$Hom_{\varinjlim_I\Cc_\bullet}(-,M)$$ and similarly, an object in $Fun(I,\vect)$ via $$Hom_{\varinjlim_I\Cc_\bullet}(M,-).$$  On the other hand, taking $\Tc=\varinjlim_I\Cc_\bullet$ and the functor to be $Id$, we see that weak sections of $\Cc_\bullet$ yield functors from $I$ to $\varinjlim_I\Cc_\bullet$.

We will next focus on trying to understand $\varinjlim_I\Cc_\bullet$ in a particular setting of direct relevance to cyclic homology.

\section{Algebras and monoidal categories}

Let us briefly recall the definition of a cyclic object associated to an algebra.  More precisely, we mean an associative algebra with $1$.  We observe that if $A$ is an algebra then this structure is encoded into a functor $A_\bullet:\Delta\to \vect$ where $\Delta$ is the simplex category of finite linearly ordered sets, with non-decreasing maps (made small).  More precisely, let $n\in \Delta$ denote the set consisting of $0,1,\cdots,n$ and set $A_n=A^{\ot(n+1)}$.  For $f\in Hom_\Delta(i,j)$ set \begin{align*}f_*(a_0\ot\cdots\ot a_i)&=\ot_{s=0}^j f^s_*(a_0\ot\cdots\ot a_i)\\f^s_*(a_0\ot\cdots\ot a_i)&=\prod_{l\in f^{-1}(s)}a_l\end{align*} where the induced ordering on $f^{-1}(s)$ is used to define the product, and $\prod_\emptyset=1$.

The above $A_\bullet$ extends to the Connes' cyclic category $\Lambda$.  One then applies cyclic duality \cite{Lod}, i.e., the identification of $\Lambda^{op}$ with $\Lambda$ by essentially exchanging the roles of faces and degeneracies to obtain $$A_\bullet:\Lambda^{op}\to\vect.$$ The latter is the classical cyclic object associated to the unital associative algebra $A$.   One obtains the simplicial object    by considering the restriction of the digram $A_\bullet$ to the simplex category $\Delta^{op}$.  The cyclic and Hochschild cohomologies of $A$ are then expressed as an $Ext$ in the abelian categories of cyclic and simplicial objects respectively \cite{ac83b}; see also \cite{kassel} for the mixed complexes version.

\subsection{Monoidal categories}\label{moncat}
Roughly speaking, a monoidal category $\Cc$ \cite{braidcat} is a category equipped with a bifunctor  $\ot:\Cc\times\Cc\to\Cc$ written as $(X,Y)\mapsto X\ot Y$ and an associativity constraint $\alpha_{X,Y,Z}:X\ot (Y\ot Z)\simeq (X\ot Y)\ot Z$ that satisfies a coherence (pentagon) axiom.  There is also a monoidal unit $1\in\Cc$ with natural requirements.  We will not be precise here since we will reinterpret this structure below.

For the sake of clarity of  exposition we begin with a construction that will be imitated below, but will not itself be immediately useful to us.   The definitions are themselves  essentially a copy of those of $A_\bullet$ above.  We define a strong $\Delta$-indexed diagram in $\Cc at$ by setting $\Cc_n=\Cc^{\times(n+1)}$.  To every $f\in Hom_\Delta(i,j)$ we assign $$f_*(X_0,\cdots, X_i)_s=\bigotimes_{l\in f^{-1}(s)}X_l,\quad\text{with}\quad\bigotimes_\emptyset=1.$$ Note that the associator eliminates the ambiguity of bracketing and defines the required structure of $\iota_{f,g}:f_*g_*\simeq (fg)_*$.  

\begin{remark}\label{multH}
There are settings, such as that of multiplier Hopf algebras \cite{multiH} when the $\iota_{f,g}$ above are not isomorphisms, so that we do not really have a monoidal category.  We nevertheless do have a weak diagram as above.
\end{remark}

Let $A$ be a unital associative algebra in $\Cc$.  The definition is exactly what one would expect, namely a monoid object in $\Cc$.  Note that $(A,\cdots, A)\in \Cc^{\times(n+1)}$ together with $f_*(A,\cdots, A)_s=\otimes_{l\in f^{-1}(s)}A\to A$ given by multiplication when $f^{-1}(s)\neq\emptyset$ and the unit, when it is empty, gives a \emph{weak} section of $\Cc_\bullet$ as defined above. 

Everything extends to diagrams and their sections over $\Lambda$.  One then applies cyclic duality to obtain a diagram and its weak section over $\Lambda^{op}$.  Recall that the input to this procedure is an algebra $A$ in a monoidal category $\Cc$.  We will take this no further here, though in principle we see that given the pair $(\Cc,A)$ we obtain a cyclic object in $\varinjlim_{\Lambda^{op}}\Cc_\bullet$.  The issue is that we don't know how to compute the limit category.  We begin to address this deficiency in the next section, by changing the limit.

\subsection{A variation on the notion of  a monoidal category}

To address the difficulty pointed out in the previous section, we need \cite{shcatchern} to replace  the category $\Cc^{\times k}$ with a more sophisticated category $\Cc^{\boxtimes k}$.  The definition of the latter is more subtle than the former, but it is worth it for us for two reasons: it is actually not that sophisticated in our examples and it makes the limit category fairly computable.

Let $\Cc=\hmod$, i.e., as an abelian category it is given as the category of left $H$-modules for some algebra $H$.  Then $\Cc^{\boxtimes 2}=\hhmod$, i.e., the category of modules with two commuting left $H$-actions; similarly $\Cc^{\boxtimes n}$ is defined as the category of modules with $n$  commuting left $H$-actions.  We have a functor $\Cc^{\times n}\to \Cc^{\boxtimes n}$  given in our example by $(X_1,\cdots, X_n)\mapsto X_1\ot_k\cdots\ot_k X_n$ where $k$ is the field over which $H$ is an algebra.  We will assume that the monoidal product on $\Cc$ extends to an appropriately associative functor $$\Delta^*:\Cc^{\boxtimes 2}\to\Cc,$$ i.e., the associator extends to $\Delta^*(Id\boxtimes\Delta^*)\simeq\Delta^*(\Delta^*\boxtimes Id)$ that satisfies the pentagon axiom.  The notation is motivated by the case of $H$ a Hopf algebra (or more generally, a bialgebra) where the product is given by the pullback along the algebra map $\Delta:H\to H\ot H$ and   its right adjoint $\Delta_*:\Cc\to \Cc^{\boxtimes 2}$, to be defined shortly. 

\begin{remark}
Going back to Remark \ref{multH} we point out that in the case of multiplier Hopf algebras the appropriate ``monoidal structure" on $\Cc$ actually depends on the desired type of the Hopf-cyclic theory.  One of the appropriate structures on $\Cc$ is monoidal and extends, the other is either weakly monoidal as in the previous remark, or strongly monoidal in the sense of this section. Thus, a weakly monoidal structure can extend to a strongly monoidal one.  We further remark that in the multiplier case the extension has a left, not right, adjoint.
\end{remark} 

We can repeat the constructions of Section \ref{moncat} to obtain a $\Lambda^{op}$-indexed diagram $\Cc^{\boxtimes(\bullet+1)}$ with a diagram morphism from the $\Cc^{\times(\bullet+1)}$ constructed therein.  Thus, the new diagram inherits a weak section (arising from $A\in\Cc$) from the old one.  The problem now reduces to the computation of $$HC(\Cc)=\varinjlim_{\Lambda^{op}}\Cc^{\boxtimes(\bullet+1)}.$$ To solve it we will assume the existence of right adjoints. This is a reasonable assumption, though sometimes we would need the left adjoints instead, the discussion would be similar.

It will be useful to us later to also consider the restriction of the diagram  $\Cc^{\boxtimes(\bullet+1)}$ to the index category $\Delta^{op}$. Let us denote the resulting colimit as follows: $$HH(\Cc)=\varinjlim_{\Delta^{op}}\Cc^{\boxtimes(\bullet+1)}.$$  Observe that the  colimit over the smaller diagram results in a category that has an extra structure when compared to the original colimit, namely an $\varsigma\in Aut(Id_{HH(\Cc)})$ that controls stability (see Definition \ref{st11}) in the sense that an object $X$ of $HH(\Cc)$ is stable if $\varsigma_X=Id$. 

\section{Towards  understanding the limit category}\label{ulc}
Let us assume that the monoidal product of $\Cc$, originally given as a bifunctor $\ot:\Cc\times\Cc\to \Cc$, not only extends to a functor $\Delta^*:\Cc\boxtimes\Cc\to\Cc$, but that the latter has a right adjoint $\Delta_*$.  Under the additional assumption that we are looking for the direct limit (colimit) defining $HC(\Cc)$ in $\Cc at_*$, i.e., categories with functors possessing right adjoints, we have that \begin{equation}\label{invlim}HC(\Cc)=\varprojlim_\Lambda\Cc_*^{\boxtimes(\bullet+1)}
\end{equation} where the diagram $\Cc_*^{\boxtimes(\bullet+1)}$ is obtained from $\Cc^{\boxtimes(\bullet+1)}$ by replacing all functors with their right adjoints.  Note that the inverse limit (limit) is computed in $\Cc at^*$.

We observe that the limit \eqref{invlim} can be explicitly described as the category of strong sections of $\Cc_*^{\boxtimes(\bullet+1)}$.  This holds in general, but in the case of our indexing category $\Lambda$ we have an even simpler description.  In the following we use the faces $\delta$'s and degeneracies $s$'s notation.  For example $\delta_0=\Delta_*$ and $\delta_1=\sigma\Delta_*$  for $\delta_i:\Cc\to\Cc^{\boxtimes 2}$.  Note that $\sigma:\Cc^{\boxtimes 2}\to \Cc^{\boxtimes 2}$ just interchanges the two $H$-actions.

A strong section, in the case of the Connes' cyclic category as index, reduces to an $M\in\Cc$ equipped  with an isomorphism $\tau$ in $\Cc^{\boxtimes 2}$ subject to conditions, namely: 
\begin{align}\label{cartsec}
\begin{split}
\tau: \delta_0(M)&\simeq\delta_1(M) \quad\text{structure in $\Cc^{\boxtimes 2}$,} \\  
s_0(\tau)&=Id_M \quad\text{unit condition in $\Cc$,}\\ 
\delta_2(\tau)\delta_0(\tau)&=\delta_1(\tau) \quad\text{associativity condition in $\Cc^{\boxtimes 3}$},\\
s_1(\tau)&=Id_M \quad\text{stability condition in $\Cc$}.
\end{split}
\end{align}

The last, stability condition, is omitted for the description of $HH(\Cc)$, instead, it equips the latter with a $$\varsigma=s_1(\tau)\in Aut(Id_{HH(\Cc)}).$$  In light of this description we see that $HH(\Cc)$ and $HC(\Cc)$ consist of objects of $\Cc$ with extra structure with the distinction between $HH(\Cc)$ and $HC(\Cc)$ being an extra property on the latter.  This is very much like the definition of the center of a monoidal category $\Zc(\Cc)$; we will come back to this in Section \ref{sec:bi}.

\begin{remark}
The  same arguments apply if we assume that the monoidal product $\ot:\Cc\times\Cc\to \Cc$ extends to $\Delta^!:\Cc\boxtimes\Cc\to\Cc$ that has a left adjoint $\Delta_!$.  The difference is illustrated in the subsequent section.
\end{remark}

\subsection{Categorical traces}\label{traces}
We observe that structure/properties \eqref{cartsec} imply that, for $M\in\Cc$ that satisfies them, the functor $$T(-)=Hom_\Cc(-,M)$$ is a (contravariant) trace  from $\Cc$ to $\vect$ \cite{trace1, trace2, trace3, trace4, hks}.   Namely, we can define natural transformations $\iota_{X,Y}:T(X\ot Y)\simeq T(Y\ot X)$, for $X,Y\in\Cc$ via \cite{shcatchern}:  \begin{align*}
Hom_\Cc(X\ot Y,M)&=Hom_\Cc(\Delta^*(X\boxtimes Y),M)=Hom_{\Cc^{\boxtimes 2}}(X\boxtimes Y,\Delta_*M)\\&\stackrel{\tau}{\to}Hom_{\Cc^{\boxtimes 2}}(X\boxtimes Y,\sigma\Delta_*M)=Hom_\Cc(\Delta^*(Y\boxtimes X),M)\\&=Hom_\Cc(Y\ot X,M).
\end{align*} These isomorphisms satisfy: $\iota_{1,X}=Id=\iota_{X,1}$, and for $X,Y,Z\in\Cc$ we have $\iota_{X\ot Y,Z}=\iota_{Y,Z\ot X}\iota_{X,Y\ot Z}$.  Note that $\iota_{1,X}=Id$ is automatic and $\iota_{X,1}=Id$ will be relaxed in Section \ref{sigmacyclic}. 
\begin{theorem}[\cite{hks}]
Any $T:\Cc\to\vect$ with this structure makes $T(A^{\ot(\bullet+1)})$, for $A\in\Cc$ an algebra, into a cocyclic object.  Similarly, for $C\in\Cc$ a coalgebra, we have that  $T(C^{\ot(\bullet+1)})$ is a cyclic object. 
\end{theorem}

Here, we were operating under the assumption that we are using the right adjoint to the product $\Delta^*$.  Similarly, we see that $L(-)=Hom_\Cc(N,-)$ is a (covariant) trace from $\Cc$ to $\vect$  if the left adjoint $\Delta_!$ to the product $\Delta^!$ is used instead.   The trace $L$ produces (co)-cyclic objects from (co)-algebras in $\Cc$.

\subsection{Monadicity} 

Let $U: HC(\Cc)\to \Cc$ denote the functor that forgets the extra structure.  Under our assumptions above  $U$ has a left adjoint $F$; the pair $(F,U)$ is monadic.  Namely, there exists a monad on $\Cc$, i.e., an algebra in $Fun(\Cc,\Cc)$, that we will denote $\Ac$, such that $HC(\Cc)$ consists of modules in $\Cc$ over $\Ac$.  Normally, monadicity is decided based on some general theorems, while the candidate monad is simply $UF$, but in our case, we are interested in an explicit description \cite{shcatchern} that we provide below.  Note that if we consider $(\Delta_!,\Delta^!)$ instead of $(\Delta^*,\Delta_*)$ then the pair $(U,F)$ is comonadic.

Observe that the only extra structure on $M\in \Cc$ that turns it into an object in $HC(\Cc)$ is the isomorphism $\tau:\Delta_*M\to\sigma\Delta_*$, while the rest are properties of $\tau$.   It is immediate that if we set $$\Ac'=\Delta^*\sigma\Delta_*$$ then $\tau$ is encoded into $act:\Ac'(M)\to M$ via adjunction.  We now need to encode the properties of $\tau$ into those of $act$.  A reasonable guess is that $\Ac'$ is a monad such that modules over it satisfy the unit and associativity conditions and furthermore, $\Ac'$ has a central element $\varsigma$ such that $s_1(\tau)$ is given by its action.  Thus, $\Ac$ should be $\Ac'/(\varsigma-1)$; this is not always true. 

The problem is not the unit and stability, those are immediate from the definition of $\Ac'$.  The problem is associativity, as in general the product on $\Ac'$ is only partially defined.

\subsection{The role of property A}\label{basechange}
If $\Cc$ is a monoidal category with the extended product that admits a right adjoint, consider the associator: $\Delta^*(Id\boxtimes\Delta^*)\simeq \Delta^*(\Delta^*\boxtimes Id)$.  Via adjunction it yields a map \begin{equation}\label{propA}(Id\boxtimes\Delta^*)(\Delta_*\boxtimes Id)\to\Delta_*\Delta^*.\end{equation}  
\begin{definition}
Roughly speaking, we say \cite{shcatchern} that $\Cc$ has property A if \eqref{propA} is an isomorphism.  In \cite{shcatchern} there are other technical conditions listed as well, but this is the key one.
\end{definition}

Recall that $\Cc$ is rigid if every object has both right and left duals, i.e., we have objects $X^*$ and ${}^*X$ such that the functor pairs $(-\ot X,-\ot X^*)$ and $(X\ot -,{}^*X\ot -)$ are adjoint pairs;  also see Remark \ref{rigbi}.

Having property A is strongly related to being rigid, in the following sense (rigid does imply it).  It is not hard to show that $\hmod$ has property A if $H$ is a Hopf algebra, but need not have it if $H$ is only a bialgebra.  The former yields a rigid category only if one considers finite dimensional $H$-modules, yet the full module category has property A as well; the bialgebra yields  a category of modules  that does not have property A even if one restricts to finite dimensional modules.

Let us consider the example of $G$, a finite monoid, i.e., a finite group without inverses.  Associate to it a monoidal category $\Cc=\vect_G$, i.e., the category of $G$-graded vector spaces with the product: \begin{equation}\label{example}(V\ot W)_g=\bigoplus_{xy=g}V_x\ot W_y.\end{equation}  Then \eqref{propA} for $\Cc$ is equivalent to the commutative diagram that expresses the associativity of the binary operation  on $G$ being Cartesian.  The latter is equivalent to $G$ being a group.  More generally, it should be a simple exercise to show that $\hmod$ has property A for a bialgebra $H$ if and only if $H$ has an invertible antipode (so that it is, in particular,  a Hopf algebra).

It turns out that if $\Cc$ has property A, then $\Ac'$ is our monad: \begin{theorem}[\cite{shcatchern}]   If $\Cc$ has property A then
 $$HH(\Cc)=\Ac'_\Cc\text{-mod}.$$ 
\end{theorem}  If $\Cc$ does not have property A, then $\Ac'$ merely generates the required monad subject to certain relations (work-in-progress of the second author).  We will see an example of this phenomenon in \eqref{mongrpd}.  

The difference between $\Delta^*\sigma\Delta_*$ itself being a monad and it merely generating one is striking.  The classical descriptions of $HH(\Cc)$ and, thus, $HC(\Cc)$ (or rather what they coincide with in Hopf-like settings) relies on $\Delta^*\sigma\Delta_*$ being a monad.  Since it fails to be one, in bialgebra-like settings, attempts to copy the classical definition \cite{kay} fail as well.  One should instead consider the monad that  $\Delta^*\sigma\Delta_*$ generates; see \eqref{mongrpd}.

\begin{remark}
If the extended product admits a left adjoint instead then we talk about comonadicity, i.e., we are looking for a comonad, whose comodules in $\Cc$ yield $HC(\Cc)$.  Everything works in the same manner.
\end{remark}

\subsection{Examples}\label{Hopfex}
Let $H$ be a Hopf algebra.  Then $\Cc=\hmod$ has property A.  The extended product in this case admits both a right and a left adjoint. By examining the resulting monad and comonad we deduce:\begin{theorem}[\cite{shcatchern}] The category $HC(\Cc)$ consists of stable anti-Yetter-Drinfeld contramodules (see \cite{contra} for the definition) in the monad case.   And $HC(\Cc)$ consists of stable anti-Yetter-Drinfeld modules in the comonad case. \end{theorem} Though the categories of stable anti-Yetter-Drinfeld contramodules and stable anti-Yetter-Drinfeld modules are  related by an adjoint pair of functors \cite{shmodcontra} that are often equivalences, they are probably not the same in general. Yet we abuse notation by speaking of the category $HC(\Cc)$ without reference to the context (assumed existence of right vs left adjoints) in which the colimits were taken.  

One should not get too attached to stable elements of $HH(\Cc)$, as we will argue in Section \ref{sigmacyclic} that it is only the first approximation to the correct setting for cyclic homology; we will replace stable objects by the so-called $\varsigma$-mixed objects.

From the general discussion it follows that for $A\in\hmod$ an algebra: \begin{equation}\label{coh}C^n=Hom_H(A^{\ot(n+1)},M)\end{equation} is a cocyclic module if $M$ is a  stable anti-Yetter-Drinfeld contramodule.  On the other hand, \begin{equation}\label{homo}C_n=Hom_H(N, A^{\ot(n+1)})\end{equation} is a cyclic module if $N$ is a  stable anti-Yetter-Drinfeld module.  

\begin{remark}
Note that \eqref{coh} and \eqref{homo} do not cover all the possible Hopf-cyclic theories.  While we did address coalgebras in Section \ref{traces}, we  omitted a discussion of certain other traces, namely ones not obtained via $Hom$ functors.  They still fit into the general framework of Section \ref{nonsense}, and are still obtained from objects in $HC(\Cc)$. 
\end{remark}

Returning to the example of $\Cc=\vect_G$ with $G$ a finite monoid, let us begin by assuming that $G$ is a group.  It is well known that in the group case, $HC(\Cc)$ consists of the full subcategory of $G$-graded $G$-equivariant vector spaces, with respect to the adjoint action of $G$ on itself; the latter is $HH(\Cc)$.  The stability condition (specifying the full subcategory) is that the action of $g$ on $V_g$ is $Id$.  The description of $HH(\Cc)$ in the group case is via the representation category of a groupoid, namely the adjoint action groupoid of $G$ on $G$.  

\begin{remark}\label{same}
Above we see a common phenomenon. Namely, as does happen in the Hopf case, outside of hard to find pathological examples \cite{hal}: $HH(\vect_G)=\Zc(\vect_G)$.  On the other hand, for bialgebras this almost never happens.
\end{remark}

To demonstrate the non property A case, consider the case of a monoid $G$.  It is not hard to show (work-in-progress of the second author), that $HH(\Cc)$ is again given as the representation category of a groupoid, but it is not an action groupoid, nor will we have $HH(\Cc)$ coincide with $\Zc(\Cc)$.   More precisely, the groupoid is defined as follows.  The objects consist of $x\in G$
 as before, but the arrows are given by pairs $(x,y)$ that have source $yx$ and target $xy$.  The composition is freely generated by these arrows subject to the relations: \begin{align}\label{mongrpd}\begin{split}(1,x)&=Id\quad\text{redundant}, \\ (x,yz)(y,zx)&=(xy,z), \quad\text{for every triple $(x,y,z)$},\\  (x,1)&=\varsigma.\end{split}\end{align}
 
 This of course reduces to the action groupoid in the group case, but for a monoid, the behaviour of the groupoid \eqref{mongrpd} is very chaotic, unlike the different groupoid whose representations give the monoidal center.  
 
 The construction of the groupoid \eqref{mongrpd} by the free generation with imposed relations, is exactly the meaning behind the statement that our monad $\Ac'$ in general is not given by $\Delta^*\sigma\Delta_*$, but is rather generated by it.  Unlike the group case, where the relations allow us to replace any sequence of generators by a single generator, the monoid case produces a rich variety of new arrows from the generating set.
 
  \subsection{Biclosed categories and centers}\label{sec:bi}
  In this section we provide a perspective  that is less general than the preceding discussion, but makes the definition of the Hochschild homology category $HH(\Cc)$ as close as possible to that of the monoidal center $\Zc(\Cc)$.  This is used in \cite{ks1, ks2}  to explore the coefficients in the Hopf algebroid, weak Hopf, and quasi-Hopf setting settings.
  
  A monoidal category $\Cc$ is biclosed if for all $X\in\Cc$, the functors $-\ot X$ and $X\ot -$ have right adjoints that we denote by $X\la -$ and $-\ra X$ respectively.  Namely, $$Hom_\Cc(X\ot Y, Z)=Hom_\Cc(X, Y\la Z)=Hom_\Cc(Y,Z\ra X).$$  
  
  \begin{remark}\label{rigbi}
  When $\Cc$ is rigid, it is biclosed with $X\la Y= Y\ot X^*$ and $Y\ra X= {}^*X\ot Y$.  Being rigid is much stronger though, as $\Cc=\hmod$ is biclosed if $H$ is a bialgebra, but only rigid if  replaced by $\Cc=\hmod_{fd}$ and $H$ has an invertible antipode.
  \end{remark}
 We can consider the center of the bimodule category \cite{egno} $\Cc^{op}$ over $\Cc$ that results from this adjoint action; let us call it $\Zc(\Cc^{op})$.  More precisely,  if $M\in \Zc(\Cc^{op})$ then for every $X\in\Cc$ we have isomorphisms: $$\iota_X: X\la M\to M\ra X,$$ subject to  $\iota_1=Id$ and $\iota_{X\ot Y}=\iota_X\iota_Y$. We observe that $$HH(\Cc)=\Zc(\Cc^{op}).$$
  
  \begin{remark}
  Suppose that $\Cc$ is pivotal, i.e., we have a monoidal natural isomorphism $Id\simeq (-)^{**}$.  Then $HH(\Cc)=\Zc(\Cc)$.  
  \end{remark}
  
 Denote by $\Zc'(\Cc^{op})$ the full subcategory of $\Zc(\Cc^{op})$ that consists of objects such that the identity map $Id\in Hom_\mathcal{C}(M,M)$ is mapped to same via \begin{equation*}\label{stabilitycond}Hom_\mathcal{C}(M,M)\simeq Hom_\mathcal{C}(1,M\la M)\to Hom_\mathcal{C}(1,M\ra M)\simeq  Hom_\mathcal{C}(M,M),\end{equation*}
 where the map in the middle is post-composition with $\iota_M$ and the isomorphisms come from the definition of the adjoint action.   	This condition is stability and so $$HC(\Cc)=\Zc'(\Cc^{op}).$$

 \subsection{The category of $\varsigma$-cyclic objects}\label{sigmacyclic}
 The goal pursued so far in this text has been to associate to an algebra $A$ in a monoidal category $\Cc$ a cyclic object in a universal manner.  This goal is achieved above by constructing just such a cyclic object in what we called $HC(\Cc)$.  Note that we did not give an explicit construction, though it can be obtained from the following discussion.
 
 The modification presented in this section  is motivated by the observation that the naive definition of $HC(\Cc)$ should be replaced by one in the setting of $\infty$-categories.  Since we are not going to pursue it at this time, we provide the discussion below as a second approximation, sufficient for the Hopf-like settings, to the truly correct setting for cyclic homology. Namely, the colimit of the diagram $\Cc^{\boxtimes(\bullet+1)}$ should have been computed in the $\infty$-setting, instead of the more comprehensible setting of the usual categories.  We propose \cite{shcatchern} that the cyclic objects in $HC(\Cc)$ should be replaced with $\varsigma$-cyclic objects in $HH(\Cc)$. \begin{definition} A $\varsigma$-cyclic object is a paracyclic  object with $$\tau_n^{n+1}=\varsigma,$$ instead of $Id$ as would be in the cyclic case. As this is a worthy replacement for $HC(\Cc)$, we will denote the category of $\varsigma$-cyclic objects in $HH(\Cc)$ by  $HH(\Cc)^{S^1}.$\end{definition}  Since  $HC(\Cc)$ consists of stable objects in $HH(\Cc)$ we retain all  the cyclic objects in $HC(\Cc)$.    We can restate this from the point of view of the mixed complexes \cite{kassel} formulation of cyclic objects.  Namely, a $\varsigma$-mixed object $(M,d,h)$ in $HH(\Cc)$ is a non-positively graded cochain complex $(M,d)$ in   $HH(\Cc)$ with a specified homotopy $h$ between $\varsigma$ and $Id$, i.e., $$dh+hd=1-\varsigma.$$  Thus, a $\varsigma$-mixed complex has cohomology in $HC(\Cc)$; the homotopy ensuring this is a very important, perhaps the most important, part of the data.  We do not distinguish between $\varsigma$-cyclic and $\varsigma$-mixed objects and denote both by $HH(\Cc)^{S^1}$.
 
 Given $A\in\Cc$ we can give a quick definition of the $\varsigma$-cyclic object in $HH(\Cc)$ associated to it.  Recall that we have an adjoint pair of functors: $$F:\Cc\leftrightarrows HH(\Cc):U.$$ We note that the functor $F$ is a $\varsigma$-trace, i.e., it is a modification of the notion of trace from Section \ref{traces} obtained by replacing $\iota_{X,1}=Id$ with $\iota_{X,1}=\varsigma_{F(X)}$.  Set \begin{equation}\label{cherndefi}ch_n(A)=F(A^{\ot(n+1)})\end{equation} to obtain the desired $\varsigma$-cyclic object.  Placing $F(A^{\ot(n+1)})$ in degree $-n$ with $d=b$ and $h=B$ (see \cite{Lod}) we get the $\varsigma$-mixed object.

 \begin{remark}
 In the Hopf algebra case we obtain a strict \cite{shmixed} generalization of stable anti-Yetter-Drinfeld modules, which have been previously used as coefficients.  Namely, we get mixed  anti-Yetter-Drinfeld modules, i.e., cochain complexes of anti-Yetter-Drinfeld modules equipped with specified homotopies between $Id$ and the stability automorphism.
 \end{remark}

\subsection{On the higher Morita invariance}\label{hmor}
Classical cyclic homology of algebras is Morita invariant, i.e., it depends not so much on algebras as on their categories of modules.  What we discuss here is the categorified version of this observation, and its generalization.  More precisely, it is well known that the monoidal center  $\Zc(\Cc)$  is invariant under duality (under suitable assumptions \cite{egno}), so that it depends only on the $2$-category of $\Cc$-modules.  Drastically different looking monoidal categories can be dual to each other, for example if $H$ is a finite dimensional Hopf algebra,   then $\hmod$ and ${}_{H^*}\Mc$ are dual.  If $H=kG$ for $G$ a finite group then the latter is $G$-graded vector spaces, and the former is $G$-representations.
 
We mentioned in Section \ref{Hopfex} that $HH(\Cc)$ yields anti-Yetter-Drinfeld contramodules for $\Cc=\hmod$, where $H$ is a Hopf algebra.  This holds in the $(\Delta^*,\Delta_*)$ setting that we fix for the purposes of this section.  It is also true that $HH(\Mc^H)$ consists of anti-Yetter-Drinfeld modules, where $\Mc^H$ is the monoidal category of right $H$-comodules.  In general we have an adjoint pair of functors \cite{shmodcontra}: \begin{equation}\label{modcontramod}(-)_c:HH(\hmod)\leftrightarrows HH(\Mc^H):\widehat{(-)},\end{equation} and they are equivalences, in particular, if $H$ is finite dimensional; these functors are compatible with stability, thus, relating the cyclic homology categories as well.

The general picture is as follows: let $\Cc$ and $\Dc$ be monoidal categories, and $\Nc$ an admissible \cite{shcatchern} $(\Cc,\Dc)$-bimodule category.  The admissibility criteria is asymmetrical, an example would be the $(\hmod,\Mc^H)$-bimodule $\vect$ that results in the \eqref{modcontramod} above.  Another, crucial, example is the $(\vect,\Cc)$-bimodule $A_\Cc\text{-mod}$ for $A\in\Cc$ algebra.  It also covers the case of a monoidal functor (with a right adjoint) $\Cc\to\Dc$ that corresponds to the $(\Cc,\Dc)$-bimodule $\Dc$.  The bimodule category $\Nc$ gives a stability compatible adjoint pair of functors: \begin{equation}\label{morita2}\Nc^*: HH(\Cc)\leftrightarrows HH(\Dc):\Nc_*\end{equation}  that can also be used to redefine \eqref{cherndefi} as $$ch(A)=A_\Cc\text{-mod}^*(k)\in HH(\Cc)^{S^1}$$ where $k$ is the trivial mixed complex.  The pair $(\Nc^*,\Nc_*)$ and \eqref{morita2} can be considered as a generalization of Morita invariance that includes cases such as \eqref{modcontramod} that are not always equivalences.   

Furthermore, for $A\in\Cc$ algebra, $\Nc^*(A)\in\Dc$  algebra can be defined, and $ch(\Nc^*(A))=\Nc^*(ch(A))$, thus, we obtain a kind of Eckmann - Shapiro lemma: \begin{theorem}[\cite{shmorita}]\label{esl}For $A\in \Cc$ and $M\in HH(\Dc)^{S^1}$ we have $$HC_\Dc^\bullet(\Nc^*A,M)=HC_\Cc^\bullet(A,\Nc_* M).$$\end{theorem}    As an example, consider an algebra $A\in\hmod$ and let $M$ be a stable anti-Yetter-Drinfeld module, then $$HC^\bullet(A,\widehat{M})=HC^\bullet(A\sharp H,M).$$

\section{Rigid braided categories}\label{rigidbraided}
In this section (it and its subsection are based on a work-in-progress of the second author) we discuss the existence of the so-called annuli stacking monoidal product \cite{annstack}, as it is known in factorization homology \cite{facthom}, on $HH(\Bc)$, for $\Bc$ rigid braided \cite{braidcat}.  Recall that a monoidal category is braided if there is a braiding isomorphism $$\tau_{X,Y}:X\ot Y\simeq Y\ot X.$$  The braiding is to satisfy certain coherence conditions that ensure that braid diagrams can be used to perform calculations; thus, this structure gives rise to braid group actions on powers $X^{\ot n}$.  The braided category $\Bc$ is symmetric if $\tau^2_{X,Y}=\tau_{Y,X}\tau_{X,Y}=Id_{X\ot Y}$; the braid action then factors through the symmetric group.  On the opposite side of the spectrum, we say that $\Bc$ is non-degenerate if $\tau^2_{X,Y}=Id_{X\ot Y}$  for all $Y$, implies that $X=1$.

The annuli stacking product structure depends on the braiding, thus, if the same monoidal category is endowed with a different braiding, the product will change.  Contrast that with the observation that the underlying category $HH(\Bc)$ needs no braiding for its definition.  When $HH(\Bc)$ is identified with $\Zc(\Bc)$, as often happens for a rigid monoidal category (this does not need a braiding), the former inherits a non-degenerate braided (pair of pants product)  from the latter; this product is unsatisfactory since it is not compatible with stability, namely $\varsigma$ is not monoidal.  The annuli stacking product on the other hand is compatible with stability, but it is not braided, unless the original braiding on $\Bc$ was symmetric.

Consider the concrete case of a Hopf algebra $H$.  The category of anti-Yetter-Drinfeld modules for $H$, unlike that of Yetter-Drinfeld modules, does not have a natural monoidal structure.  It is, nevertheless, a bimodule category over the braided category of Yetter-Drinfeld modules; see Lemma \ref{lem43}, \cite{hkrs1}. It is known that often, in fact, in the presence of the original modular pair in involution, one has an identification between  the two categories.  It is even possible to describe the action of $\varsigma$ in terms of a certain ribbon element of the category of Yetter-Drinfeld modules.   This identification endows, however non-canonically, the  anti-Yetter-Drinfeld modules with a product.  Unfortunately $\varsigma$ is not monoidal with respect to this product. In fact, any ribbon element is not monoidal (it is a twist, see \eqref{twist:def}) unless the braided category is symmetric, which the category of Yetter-Drinfeld modules never is.  Thus, the identification produced by a modular pair in involution does not produce a product on the stable anti-Yetter-Drinfeld modules; also, it is rather arbitrary.
 
We can demonstrate the contrast between the two types of product on $HH(\Bc)$ in our example \eqref{example}.  We assume that $G$ is a group here.  It is true that $\vect_G$ is not braided (as $G$ need not be commutative) but the category $Rep(G)$ of representations of $G$ is (in fact it is symmetric), and they share the same Hochschild homology category, as in Section \ref{hmor}.  Recall that $\Zc(\vect_G)=HH(\vect_G)=Sh_G G^{ad}$, where the latter is $G$-graded, $G$-equivariant (with respect to the $x(-)x^{-1}$ action) vector spaces.  The product on $\Zc(\vect_G)$ is the convolution product, i.e, it is given by the same formula is in \eqref{example}.  On the other hand the correct product on $HH(\vect_G)$ is given by \begin{equation}\label{otherprod}(V\ot W)_x=V_x\ot W_x,\end{equation} which is clearly compatible with stability, i.e., $\varsigma(v_x)=xv_x$ is monoidal with respect to \eqref{otherprod}, but is a twist \eqref{twist:def} with respect to \eqref{example}.
 
We describe, roughly, a construction of the product on $HH(\Bc)$. Let $\Bc$ be rigid braided, the category $HH(\Bc)$ is given not just by a monad on $\Bc$ as was discussed in Section \ref{basechange}, but (using the braiding) by a commutative algebra $\Hc$ in $\Zc(\Bc)$.  By a generalization of the usual monoidal product on modules over a commutative ring, we obtain the desired product on $\Hc_\Cc\text{-mod}=HH(\Bc)$.  
 
 It is possible to describe the product more concretely in our case of interest (see Section \ref{taft}).  We will require an additional assumption that $\Bc$ is balanced, i.e, it is equipped with a twist $\theta\in Aut(Id_\Bc)$ satisfying:  \begin{equation}\label{twist:def}\theta_{X\ot Y}=\tau^2_{X,Y}\theta_X\ot\theta_Y.\end{equation} 
 Now, if the braiding on $\Bc$ is non-degenerate (the opposite of symmetric) then \begin{equation}\label{endo}HH(\Bc)=End(\Bc),\end{equation} as monoidal categories, in contrast with our example above where $\Bc$ was symmetric.
We can  extend \eqref{endo} in a way that elucidates $HC(\Bc)$.    Namely, under the equivalence  in \eqref{endo} we have the correspondence $$\varsigma=\theta(-)\theta^{-1}.$$  Note that $\Bc$ may have more than one twist (usually many as they form a torsor over $Aut^\ot(Id_\Bc)$), yet $\varsigma$ is $\varsigma$.  The role of different twists will become evident below (see Remark \ref{twistscoh}).
 
 It may seem, especially in light of \eqref{endo}, that this canonical product on $HH(\Bc)$ is not very interesting.  This is not the case; recall that this product is compatible with stability and, thus, restricts to $HC(\Bc)$.  The latter is more interesting.  Furthermore, the product plays a crucial role in the next section. 
 
 \subsection{Relative coefficients and relative centers}
 Consider a Hopf algebra $H$ in a rigid braided category $\Bc$ \cite{majidbraid}.  It is immediate that $$\Cc=H_\Bc\text{-mod},$$ the category of  modules over $H$ in $\Bc$, is a monoidal category.  We can, thus, apply the machinery considered in this paper to $\Cc$, and obtain $HH(\Cc)$ and $HC(\Cc)$, but it would not remember $\Bc$.  In fact, $\Cc$ can often be expressed as modules of a certain ``product" Hopf algebra $T$ (see \eqref{tafteq} below), and the machinery does not distinguish between  $H_\Bc\text{-mod}$ and $T\text{-mod}$.
 
 What  one wants instead is a relative (with respect to $\Bc$) version of $HH(\Cc)$ and $HC(\Cc)$.   In fact, one may want a relative version of the more classical $\Zc(\Cc)$ as well (this is considered in \cite{laug}).  The latter problem is actually simpler, so we will briefly address it as well.  It is tempting, as was attempted in \cite{kp}, to try and generalize the old constructions to yield new (anti)-Yetter-Drinfeld modules, now starting with a braided category other than the $\vect$ of before.   This works for Yetter-Drinfeld modules, and, in principle,  for anti-Yetter-Drinfeld modules in a symmetric category, but it cannot work in general. One needs ``coefficients" for this to work. 
 
We begin by observing that if $H$ is a Hopf algebra in $\Bc$ then we have a braided functor $$\Bc\to\Zc(\Cc).$$  This induces a monoidal, with respect to the annuli stacking product,  functor $$HH(\Bc)\to HH(\Zc(\Cc)).$$  The latter monoidal category acts on the right on the  category $HH(\Cc)$.  We see that $HH(\Cc)$ is a right $HH(\Bc)$-module; the same applies to $HC(-)$.  In particular, if $\Mc$ is a left $HH(\Bc)$-module (with a compatible $\varsigma$-action) then one may consider the category \begin{equation}\label{relcoh}HH(\Cc)\boxtimes_{HH(\Bc)}\Mc\end{equation} or its cyclic version: $HC(\Cc)\boxtimes_{HC(\Bc)}\Mc^\varsigma.$  This right module structure of $HH(\Cc)$ over $HH(\Bc)$ is the ``memory" of $\Bc$.  As we see in \eqref{localrem}, it localizes $HC(\Cc)$.

When considering the center version of the above we note that the braided functor $\Bc\to\Zc(\Cc)$ also yields a monoidal, with respect to the annuli stacking product (not the usual product on centers), functor $\Zc(\Bc)\to \Zc(\Zc(\Cc))$, and the latter acts on the right on $\Zc(\Cc)$.  In this case $\Zc(\Bc)$ (with  the annuli stacking product) has a canonical left module category $\Bc$, so that the relative center can be defined as \begin{equation}\label{relcent}\Zc_\Bc(\Cc)=\Zc(\Cc)\boxtimes_{\Zc(\Bc)}\Bc.\end{equation}

Going back to our Hochschild homology case, the action of $HH(\Bc)$ on $\Bc$ is equivalent to specifying a twist $\theta$.  If $\Bc$ is symmetric then the trivial twist is available and so the relative Hochschild and cyclic homology categories can be defined as in   \eqref{relcent}.  In the non-symmetric case, there is no canonical twist (or none may be available).  However, with a specified twist in hand we can define the relative Hochschild and cyclic homology categories: \begin{equation}\label{twistcoh}HH_{(\Bc,\theta)}(\Cc)=HH(\Cc)\boxtimes_{HH(\Bc)}\Bc_\theta\quad\text{and}\quad HC_{(\Bc,\theta)}(\Cc)=HC(\Cc)\boxtimes_{HC(\Bc)}\Bc_\theta^\theta.\end{equation} 
\begin{remark}
The definition \eqref{relcoh} is rather abstract, but it is possible to describe these categories monadically.  The resulting definition \eqref{twistcoh}  of relative anti-Yetter-Drinfeld modules (dependent on the $\theta$ chosen) looks very similar to the old one.  More precisely, they are described as objects in $\Bc$ that are modules over $H$ and ${}^*H$ such that the actions satisfy a compatibility condition  that involves $\theta$.  The definition of $\varsigma$ also involves $\theta$.  Unlike the $\vect$ case where anti-Yetter-Drinfeld modules are modules over the twisted Drinfeld double, here we have that they are modules over a suitable monad on $\Bc$, not an algebra in $\Bc$.
\end{remark}

\subsection{Examples}\label{taft}
Let $G$ be an abelian group and $\chi$ a bicharacter such that $\omega(x,y)=\chi(x,y)\chi(y,x)$ is non-degenerate.  Let $\Bc=(\vect_G,\chi)$, i.e., $$\tau: v_x\ot w_y\mapsto \chi(x,y)w_y\ot v_x$$ and $\theta(g):=\chi(g,g)$ its canonical  ribbon element (special twist). Let $$I=\theta(G)\subset k^\times\quad\text{and}\quad n_i=\#\{x|\theta(x)=i\}$$ then \begin{equation}\label{eq:gab}HH(\Bc)\simeq M_{|G|}(\vect)\quad\text{and}\quad HC(\Bc)=\bigoplus_{i\in I}M_{n_i}(\vect),\end{equation} where $M_n(\vect)$ is the monoidal category of  matrices with vector space entries.  Any other twist is of the form $\theta_x=\theta\omega(x,-)$ for $x\in G$.
 The twists $\theta_x$ and $\theta_y$ yield the same theory  if $\theta(x)=\theta(y)$.  Thus, for an $H\in\Bc$ there is a unique category of ``relative" anti-Yetter-Drinfeld modules, but $\# I$ categories of ``relative" stable anti-Yetter-Drinfeld modules in general.

\begin{remark}\label{twistscoh}
In \eqref{eq:gab} we see that $HC(\Bc)$ has irreducible module categories $\vect^{n_i}$ indexed by $i\in I$.  These correspond to twists $\theta_x$ with $\theta(x)=i$.
\end{remark}
 
We can go further in a subexample of the above by considering Taft algebras \cite{taft}. We recall the relevant definitions: suppose that $p$ is an odd prime and $\xi$ is a primitive $p$th root of $1$.  Let $\Bc=(\vect_{\mathbb{Z}/p},\xi^{ij})$, and enumerate the twists by $\theta_t=\xi^{i^2+ti}$.  Define  $H=k[x]/x^p$ with $deg(x)=1$,  $\Delta(x)=1\ot x+ x \ot 1$, $S(x)=-x$, and $\epsilon(x)=0$.   Note that \begin{equation}\label{tafteq}H_\Bc\text{-mod}=T_p(\xi)\text{-mod},\end{equation} where  $T_p(\xi)$ is the Taft Hopf algebra.
 Let $q=\xi^{-1/2}$ then  independently of $t$: $$HH_{(\Bc,\theta_t)}(T_p(\xi)\text{-mod})\simeq u_q(sl_2)\text{-mod},$$ where $u_q(sl_2)$ denotes the small quantum group, while $$HC_{(\Bc,\theta_0)}(T_p(\xi)\text{-mod})\simeq \vect.$$
  For $t\neq 0$,  $HC_{(\Bc,\theta_t)}(T_p(\xi)\text{-mod})$ depends only on $|t|$ but is otherwise non-trivial and sits over a distinct (for $|t|\neq |t'|$) point of $Spec \mathcal{Z}(u_q(sl_2))$.

In light of the above,  the question of what is a Hopf-cyclic theory for a Hopf algebra in a braided category has two answers.  First, balanced (not braided) categories yield Hopf-cyclic theories.  Different twist = Different theory.  Second, remembering $\Bc$ ``localizes" $HC(\Cc)$.  In the Taft algebra example (which is just a usual Hopf algebra) we have that any stable anti-Yetter-Drinfeld module (in the classical sense)  breaks down $$M=M_0\oplus\bigoplus_{i=1}^{(p-1)/2}\left(M_i^1\oplus M^2_i\right)$$ with $M_0$ in a category equivalent to $\vect$ and $M^j_i$ of type $i$.  The index  $0,\cdots,(p-1)/2$ runs over the spectrum of $\mathcal{Z}(u_q(sl_2))$; the type $i\neq 0$ categories are not $\vect$.

Returning to the non-degenerate case of $\Bc=(\vect_G,\chi)$ with $\Cc=H_\Bc\text{-mod}$, there exists a category $X$ with $\varsigma\in Aut(Id_X)$ such that \begin{equation}\label{localrem}HC(\Cc)=\bigoplus_{i\in I}X^{i\varsigma}\boxtimes\vect^{n_i}.\end{equation}

\bibliographystyle{amsalpha}

\end{document}